\documentclass[leqno,11pt]{amsart}
\usepackage{amsmath,amscd}
\usepackage{epsfig,graphics}
\usepackage{amssymb,latexsym}
\usepackage{mathrsfs}
\newtheorem{thm}{\bf Theorem}[section]
\newtheorem{df}[thm]{\bf Definition}
\newtheorem{prop}[thm]{\bf Proposition}
\newtheorem{cor}[thm]{\bf Corollary}
\newtheorem{lem}[thm]{\bf Lemma}
\newtheorem{rem}[thm]{\bf Remark}
\newtheorem{ex}[thm]{\bf Example}
\newcommand{\A}{\mathcal{A}}
\newcommand{\B}{\mathcal{B}}
\newcommand{\W}{\mathcal{W}}
\newcommand{\cP}{\mathscr{P}}
\newcommand{\pf}{\noindent{\bfseries Proof. }}

\numberwithin{equation}{section}

\begin{document}
\title[]
{A combinatorial proof of a Weyl type formula for hook Schur
polynomials}
\author{JAE-HOON KWON}
\address{Department of Mathematics \\ University of Seoul \\ 90
Cheonnong-dong, Dongdaemun-gu \\ Seoul 130-743, Korea }
\email{jhkwon@uos.ac.kr }

\thanks{This research was supported by KRF Grant $\sharp$2005-070-C00004 } 

\maketitle

\begin{abstract}
In this paper, we present a simple combinatorial proof of a Weyl
type formula for hook Schur polynomials, which has been obtained by
using a Kostant type cohomology formula for $\frak{gl}_{m|n}$. In
general, we can obtain in a combinatorial way a Weyl type formula
for various highest weight representations of a Lie superalgebra,
which together with a general linear algebra forms a Howe dual pair.
\end{abstract}

\section{Introduction}

The notion of hook Schur polynomial was introduced by Berele and
Regev \cite{BR}, as characters of complex irreducible tensor
representations of the general linear Lie superalgebra
$\frak{gl}_{m|n}$.  Recently, in \cite{CZ} Cheng and Zhang proved a
Kostant type cohomology formula for $\frak{gl}_{m|n}$ associated to
its irreducible tensor representations to compute the corresponding
generalized Kazhdan-Lusztig polynomials (cf.\cite{S}), which also
implies a Weyl type formula for irreducible tensor representations
by Euler-Poincar\'{e} principle (cf.\cite{GL}). This Weyl type
formula, which is given as an alternating sum of characters of Kac
modules, is closely related with a general approach to the study of
the complex irreducible finite dimensional representations of
$\frak{gl}_{m|n}$ (cf.\cite{Br,Kac78,S}).

In this paper, we introduce a new combinatorial proof of the Weyl
type formula for hook Schur polynomials  obtained in \cite{CZ}. Our
proof, which was originally motivated by \cite{CZ-2}, is  simple and
natural in the sense that we use only the classical Weyl formula and
the Cauchy identity for Schur polynomials.  In fact, we prove a Weyl
type formula for a more general class of functions which arise
naturally as characters of quasi-finite irreducible representations,
not necessarily finite dimensional, of various Lie (super)algebras
\cite{CW03,EHW,Fr,KacR2}, and were also introduced in a
combinatorial way in \cite{Kwon}.

Suppose that $\A$ and $\B$ are $\mathbb{Z}_2$-graded sets at most
countable, and $\lambda$ is a generalized partition of length $d$.
Let
\begin{equation*}
S_{\lambda}^{\A/\B}=\sum_{\mu,\nu}S_{\mu}({\bf x}_{\A}) S_{\nu}({\bf
x}_{\B}^{-1}),
\end{equation*}
where $S_{\mu}({\bf x}_{\A})$ and $S_{\nu}({\bf x}_{\B}^{-1})$ are
super Schur functions (or super symmetric functions) in the
variables ${\bf x}_{\A}=\{\,x_a\,|\,a\in\A\,\}$ and ${\bf
x}_{\B}^{-1}=\{\,x_b^{-1}\,|\,b\in\B\,\}$ corresponding to skew
shapes $\mu=\left(\lambda+(p^d)\right)/\eta$ and $\nu=(p^d)/\eta$
for some $p\geq 0$ and $\eta$ (see Definition \ref{SAB}). Then the
main result (Theorem \ref{main result}) is
$$S_{\lambda}^{\A/\B}=
\frac{\sum_{w\in \W} (-1)^{\ell(w)} S_{\lambda^{w,+}}({\bf
x}_{\A})S_{\lambda^{w,-}}({\bf x}_{\B}^{-1})}{\Delta_{\A/\B}},$$
where $\Delta_{\A/\B}=\prod_{|a|=|b|}(1-x_ax_b^{-1})\prod_{|a|\neq
|b|}(1+x_ax_b^{-1})^{-1}$, $\W$ is a set of right coset
representatives of an affine Weyl group of type $A_{\infty}$ with
respect to a maximal parabolic subgroup, and $\lambda^{w,\pm}$ are
defined under an action of $w\in \W$ on $\lambda$. We also give
alternative proofs of a Cauchy identity of $S^{\A/\B}_{\lambda}$
paired with rational Schur polynomials and a Jacobi-Trudi identity
for $S^{\A/\B}_{\lambda}$ (cf.\cite{Kwon}). Some of the arguments
might be stated or understood more easily in the context of
representation theory, but we  give self-contained combinatorial
proofs which do not depend on it.

Now a Weyl type formula for hook Schur polynomials or irreducible
tensor representations of $\frak{gl}_{m|n}$ (Theorem \ref{Weyl for
hook Schur}) follows as a byproduct, up to a multiplication of a
monomial, when $\A$ and $\B$ are finite sets (say $|\A|=n$ and
$|\B|=m$) of even and odd degree, respectively.  This recovers in a
purely combinatorial way the character formula given in \cite{CZ}.
We also give another proof of the factorization property of hook
Schur polynomials (cf.\cite{BR,rem}).

In general, we can obtain Weyl type character formulas for other
irreducible highest weight representations of a Lie (super)algebra,
whenever it forms a  Howe dual pair with a general linear algebra
(cf. \cite{CW01,CW03,Fr,H,KacR2}), since the characters of the
associated representations are given by $S_{\lambda}^{\A/\B}$ under
suitable choices of $\A$ and $\B$ \cite{Kwon}. We discuss in detail
one more example in representation theory when both $\A$ and $\B$
are finite sets of even degree (say $|\A|=n$ and $|\B|=m$). We
deduce from Howe duality \cite{H,KV} that the corresponding
$S^{\A/\B}_{\lambda}$, up to a multiplication of a monomial, is a
character of an infinite dimensional representation of
$\frak{gl}_{m+n}$, which is of particular importance in the study of
unitary highest weight representations of the Lie group $U(m,n)$
(cf.\cite{EHW}). In this case, we obtain a Weyl type formula given
as an alternating sum of characters of generalized Verma modules
(Theorem \ref{Weyl for unitary}), which recovers the Enright's
character formula \cite{En} with a different parametrization of
highest weights for generalized Verma modules, and also an analogue
of the Jacobi-Trudi formula for these infinite dimensional
representations.

Finally, we would like to mention that the similarity of character
formulas for $\frak{gl}_{m|n}$ in \cite{CZ} and $\frak{gl}_{m+n}$ in
\cite{En} was already observed, and a more direct connection between
the Grothendieck groups of module categories of $\frak{gl}_{m|n}$
and $\frak{gl}_{m+n}$ has been established in \cite{CWZ} recently.

The paper is organized as follows. In Section 2, we recall some
basic terminologies. In Section 3, we derive a Weyl type formula, a
Cauchy type identity, and a Jacobi-Trudi formula for
$S^{\A/\B}_{\lambda}$. We also discuss a factorization property for
$S^{\A/\B}_{\lambda}$ when $\A$ is a finite set of even degree,
which provides another proof of the factorization of hook Schur
polynomials. In Section 4, we discuss applications to irreducible
tensor representations of $\frak{gl}_{m|n}$, and infinite
dimensional representations of $\frak{gl}_{m+n}$.\vskip 3mm

{\bf Acknowledgement} Part of this work was done during the author's
visit at National Taiwan University in 2006 summer. He would like to
thank Shun-Jen Cheng for the invitation and kind explanations on his
recent works.

\section{Symmetric function}

Let us recall some  terminologies (cf.\cite{Mac}). A partition is a
non-increasing sequence of non-negative integers $\lambda =
(\lambda_k)_{k\geq 1}$ such that $\sum_{k\geq 1}\lambda_k < \infty$.
The number of non-zero parts in $\lambda$ is called the length of
$\lambda$ denoted by $\ell(\lambda)$. We also write
$\lambda=(1^{m_1},2^{m_2},\ldots)$, where $m_i$ is the number of
parts equal to $i$. We denote by $\cP$ the set of all partitions. A
partition $\lambda = (\lambda_k)_{k\geq 1}$ is identified with a
Young diagram which is a collection of nodes (or boxes) in
left-justified rows with $\lambda_k$ nodes in the $k$th row. We
denote by $\lambda'$ the conjugate of $\lambda$. For
$\lambda,\mu\in\cP$, let $\lambda+\mu=(\lambda_k+\mu_k)_{k\geq 1}$,
and if $\lambda\supset\mu$ (that is, $\lambda_k\geq \mu_k$ for all
$k$), let $\lambda/\mu$ be the skew Young diagram obtained from
$\lambda$ by removing $\mu$.

For a set $K$ which is at most countable, let $\Lambda_K$ be the
ring of symmetric functions in the variables ${\bf
x}_K=\{\,x_k\,|\,k\in K\,\}$, and $s_{\lambda}({\bf x}_K)$ the Schur
function corresponding to $\lambda\in \cP$. When $K$ is infinite,
let $\omega_K$ be the involution on $\Lambda_K$, which sends
$s_{\lambda}({\bf x}_K)$ to $s_{\lambda'}({\bf x}_K)$.

Throughout the paper, we denote by $\A=\A_0\sqcup\A_1$ a
$\mathbb{Z}_2$-graded set, which is at most countable. For $a\in
\A$, $|a|$ denotes the degree of $a$. We put $\mathbb{Z}=\{\,0,\pm
1,\pm 2,\ldots\,\}$, $\mathbb{Z}_{> 0}=\{\,1,2,\ldots\,\}$,
$\mathbb{Z}_{< 0}=\{\,-1,-2,\ldots\,\}$, $[n]=\{\,1,\ldots,n\,\}$,
and $[-n]=\{\,-1,\ldots,-n\,\}$ ($n\geq 1$), where all the elements
are assumed to be of degree $0$ (or even). Also, we define
$\A'=\{\,a'\,|\,a\in \A\,\}$ to be the set with the opposite
$\mathbb{Z}_2$-grading, that is, $|a'|\equiv |a|+1 \pmod{2}$ for
$a\in\A$.

Let ${\bf x}^{\pm 1}_{\A}=\{\,x^{\pm 1}_a\,|\,a\in\A\,\}$  be the
set of variables indexed by $\A$.  For a skew Young diagram
$\lambda/\mu$, a super Schur function corresponding to $\lambda/\mu$
is defined to be
\begin{equation}
S_{\lambda/\mu}({\bf x}_{\A})=\sum_{\substack{\nu\in \cP
\\ \mu\subset\nu\subset\lambda }}s_{\nu/\mu}({\bf x}_{\A_0})s_{\lambda'/\nu'}({\bf x}_{\A_1})
\end{equation}
(cf.\cite{BR,Mac}). For simplicity, let us often write
$S^{\A}_{\lambda/\mu}=S_{\lambda/\mu}({\bf x}_{\A})$. When $\A$ is
finite, $S^{\A}_{\lambda}$ is a hook Schur polynomial introduced by
Berele and Regev \cite{BR}. Following our notation, we may write
$S^{\A}_{\lambda}=\sum_{\mu\in
\cP}S_{\mu}^{\A_0}S_{\lambda/\mu}^{\A_1}$ for $\lambda\in \cP$, and
hence $S^{\A}_{\lambda}=s_{\lambda}({\bf x}_{\A})$ if $\A=\A_0$, and
$S^{\A}_{\lambda}=s_{\lambda'}({\bf x}_{\A})$ if $\A=\A_1$.

For a positive integer $d$, let
$\mathbb{Z}^d_+=\{\,\lambda=(\lambda_1,\cdots,\lambda_d)\in\mathbb{Z}^d\,|\,\lambda_1\geq
\cdots\geq \lambda_d\,\}$ be the set of generalized partitions of
length $d$. Put
\begin{equation}\label{lambda}
\begin{split}
\lambda^+&=(\max(\lambda_1,0),\ldots,\max(\lambda_d,0))\in \cP, \\
\lambda^-&=(\max(-\lambda_d,0),\ldots,\max(-\lambda_1,0))\in \cP,
\\
\lambda^*&=(-\lambda_d,\ldots,-\lambda_1)\in\mathbb{Z}_+^d.
\end{split}
\end{equation}
The addition on $\mathbb{Z}_+^d$ is defined in a usual way, and then
$\lambda=\lambda^+ +(\lambda^-)^*$.

 \begin{df}[\cite{Kwon}]\label{SAB}{\rm
Let $\A$ and $\B$ be $\mathbb{Z}_2$-graded sets, which are at most
countable. For $\lambda\in\mathbb{Z}_+^d$, we define
\begin{equation*}
S_{\lambda}({\bf x}_{\A};{\bf x}_{\B})=\sum_{\mu,\nu}S_{\mu}({\bf
x}_{\A}) S_{\nu}({\bf x}_{\B}^{-1}),
\end{equation*}
where $\mu$ and $\nu$ are skew Young diagrams of the form
\begin{equation*}
\mu=\left(\lambda+(p^d)\right)/\eta, \ \ \ \nu=(p^d)/\eta
\end{equation*}
for some non-negative integer $p$ and partition $\eta$ such that
$\lambda+(p^d)\in\cP$ and $\eta\subset \lambda+(p^d),(p^d)$. Let us
write $S^{\A/\B}_{\lambda}=S_{\lambda}({\bf x}_{\A};{\bf x}_{\B})$
for simplicity.}
\end{df}
\begin{rem}{\rm
(1) If $\B$ is empty, then $S^{\A/\B}_{\lambda}$ is non-zero only if
$\lambda$ is an ordinary partition, and in this case, we have
$S^{\A/\B}_{\lambda}=S^{\A}_{\lambda}$.

(2) By definition, $S^{\A/\B}_{\lambda}$ can be regarded as the
character of certain bitableaux. A combinatorics of these
bitableaux, including analogues of the Schensted insertion, the
Littlewood-Richardson rule, and the Robinson-Schensted-Knuth
correspondence, are given in \cite{Kwon}.}
\end{rem}

For $\lambda\in\mathbb{Z}_+^d$, let $s_{\lambda}({\bf x}_{[d]})$ be
the rational Schur polynomial corresponding to $\lambda$, that is,
$s_{\lambda}({\bf x}_{[d]})=(x_1\cdots
x_d)^{-p}s_{\lambda+(p^d)}({\bf x}_{[d]})$ for $p\geq 0$ such that
$\lambda+(p^d)\in\cP$. For $\mu,\nu\in\mathbb{Z}_+^d$, we have
\begin{equation}\label{LR rule for rational Schur}
s_{\mu}({\bf x}_{[d]})s_{\nu}({\bf
x}_{[d]})=\sum_{\lambda\in\mathbb{Z}_+^d}c^{\lambda}_{\mu\,\nu}s_{\lambda}({\bf
x}_{[d]}),
\end{equation}
where $c^{\lambda}_{\mu\,\nu}$ is a Littlewood-Richardson
coefficient. Note that
$c^{\lambda}_{\mu\,\nu}=c^{\lambda+((p+q)^d)}_{\mu+(p^d)\,\nu+(q^d)}$
for all  $p,q\geq 0$. Then we have another expression of
$S^{\A/\B}_{\lambda}$ as a linear combination of the products
$S_{\mu}({\bf x}_{\A})S_{\nu}({\bf x}_{\B}^{-1})$ for $\mu,\nu\in
\cP$.
\begin{prop}[\cite{Kwon}]\label{character for SAB} For  $\lambda\in\mathbb{Z}_+^d$, we have
\begin{equation*}
S_{\lambda}({\bf x}_{\A};{\bf x}_{\B})=\sum_{\substack{\mu,\nu\in
\cP \\ \ell(\mu),\ell(\nu)\leq
d}}c^{\lambda}_{\mu\,\nu^*}S_{\mu}({\bf x}_{\A})S_{\nu}({\bf
x}^{-1}_{\B}).
\end{equation*}
\end{prop}\qed

\section{Weyl type formula}

\subsection{Main result}
Put $\mathbb{Z}^{\times}=\mathbb{Z}\setminus\{0\}$. Let
$P=\bigoplus_{i\in\mathbb{Z}^{\times}}\mathbb{Z}\epsilon_i$ be the
free abelian group generated by
$\{\,\epsilon_i\,|\,i\in\mathbb{Z}^{\times}\,\}$. For
$i\in\mathbb{Z}$, let $r_i$ be the transposition on
$\mathbb{Z}^{\times}$ (hence on
$\{\,\epsilon_i\,|\,i\in\mathbb{Z}^{\times}\,\}$) given by
\begin{equation}
r_i=
\begin{cases}
(i\ i+1), & \text{if $i>0$}, \\
(i-1\ i), & \text{if $i<0$}, \\
(-1\ 1), & \text{if $i=0$}.
\end{cases}
\end{equation}
Let $W$ be the Coxeter group generated by
$\{\,r_i\,|\,i\in\mathbb{Z}\,\}$, and $\ell(w)$ denotes the length
of $w\in W$. For each subset $I\subset \mathbb{Z}$, let $W_I$ be the
subgroup of $W$ generated by $\{\,r_i\,|\,i\in I\,\}$. Let
\begin{equation}
\W=\{\,w\in W\,|\,\ell(r_iw)>\ell(w) \ \text{for $i\in
\mathbb{Z}^{\times}$}\,\},
\end{equation}
be the set of right coset representatives with respect to a maximal
parabolic subgroup $W_{\mathbb{Z}^{\times}}$ (cf.\cite{BB}).

For $\lambda \in\mathbb{Z}_+^d$ and $w\in \W$, choose sufficiently
large $p,q>0$ such that
\begin{itemize}
\item[(1)] $ -p\leq \lambda_d\leq\lambda_1\leq q$,

\item[(2)] $w\in W_{I(p,q)}$,
where $I(p,q)=\{\,k\in\mathbb{Z}\,|\,-p+1\leq k\leq q-1\,\}$.
\end{itemize}
Put $\mu=(\lambda+(p^d))'=(\mu_1,\ldots,\mu_n)$, where $n=p+q$. We
may identify $\mu$ with
$$\mu=\mu_1\epsilon_{-p}+\cdots+\mu_p\epsilon_{-1}+\mu_{p+1}\epsilon_1+\cdots+\mu_{p+q}\epsilon_q\in P.$$
Then, we define
\begin{equation}
w\circ \lambda = w(\mu+\delta_{p,q})-\delta_{p,q}-d{\bf 1}^-_p,
\end{equation}
where $\delta_{p,q}=\sum_{i\in [-p]}(q-i-1)\epsilon_{i} +\sum_{j\in
[q]}(q-j)\epsilon_j$, and ${\bf 1}^-_p=\sum_{i\in
[-p]}\epsilon_{i}$.

\begin{lem}\label{shifted action} Under the above hypothesis,
there exist unique $\sigma, \tau\in\cP$ with $\ell(\sigma)\leq p$
and $\ell(\tau)\leq q$ such that
$$w\circ \lambda= -\sum_{i\in \mathbb{Z}_{<0}}\sigma_{-i}\epsilon_{i}+
\sum_{j\in \mathbb{Z}_{>0}}\tau_{j}\epsilon_j.$$
\end{lem}
\pf First, let
$w(\delta_{p,q})-\delta_{p,q}=\sum_{i\in\mathbb{Z}^{\times}}a_i\epsilon_i
\in P$, where $a_i=0$ for $i\not\in [-p]\cup [q]$. It is not
difficult to see that
$w(\delta_{p',q'})-\delta_{p',q'}=w(\delta_{p,q})-\delta_{p,q}$ for
all $p'>p$ and $q'>q$. Since $w$ is a right-coset representative of
$W_{\mathbb{Z}^{\times}}$ in $W$, we also have $a_i\geq a_{i+1}$ and
$a_{-i-1}\geq a_{-i}$ for all $i> 0$ (cf.\cite{BB}). This implies
that
\begin{equation}
\begin{split}
& 0\geq a_{-p}\geq a_{-p+1}\geq \ldots \geq a_{-1}, \\
& a_1\geq a_2\geq \ldots \geq a_q \geq 0.
\end{split}
\end{equation}
Next, if we put $w(\mu)=\sum_{i\in [-p]\cup [q]}b_i\epsilon_i$, then
we have
\begin{equation}
d\geq b_{-p}\geq \cdots\geq b_{-1}\geq 0, \ \ b_1\geq \cdots\geq
b_q\geq 0.
\end{equation}
Note that $w(\mu)-d{\bf 1}^-_p$ does not depend on the choice of
$p,q$. Hence it follows that $w\circ \lambda=\sum_{i\in [-p]\cup
[q]}c_i\epsilon_i$, where
\begin{equation}
\begin{split}
& 0\geq c_{-p}\geq c_{1-p}\geq \ldots \geq c_{-1}, \\
& c_1\geq c_2\geq \ldots \geq c_q \geq 0.
\end{split}
\end{equation}
This completes the proof. \qed\vskip 3mm

\begin{df}{\rm
For $\lambda\in\mathbb{Z}_+^d$ and $w\in \W$, we define
\begin{equation}\label{lambda w}
\lambda^{w,-}=\sigma', \ \ \ \ \lambda^{w,+}=\tau',
\end{equation}
where $\sigma,\tau\in\cP$ are given in Lemma \ref{shifted action}. }
\end{df}

\begin{rem}\label{comments on cosets}{\rm
(1) Given $w\in\W$, suppose that $w\in W_{I(p,q)}$ for some $p,q>0$.
If
$w(\delta_{p,q})-\delta_{p,q}=\sum_{i\in\mathbb{Z}^{\times}}\mu_i\epsilon_i
\in P$, then we can check that the partition $(-\mu_{-k})_{k\geq 1}$
is the conjugate of $(\mu_k)_{k\geq 1}$ (see 2.4 in \cite{BB}).
Moreover, the map sending $w$ to $(\mu_k)_{k\geq 1}$ is a one-to-one
correspondence between $\W$, the set of the minimal length right
coset representatives and $\cP$, where $\ell(w)=|\mu|=\sum_{k\geq
1}\mu_k$.

(2) One can also check that for $\lambda\in\mathbb{Z}_+^d$,
$(\lambda^{w,-}, \lambda^{w,+})=(\lambda^{w',-}, \lambda^{w',+})$ if
and only if $w=w'\in \W$.

}
\end{rem}

\vskip 3mm

Next, consider the Schur polynomials in $n$ variables. Fix $p,q>0$
such that $p+q=n$. Instead of $[n]$, let us use $[-p,q]=[-p]\cup
[q]$ as an index set for the variables. For a partition $\mu$ with
$\ell(\mu)\leq n$, we may identify $\mu$ with
$$\mu_1\epsilon_{-p}+\cdots+\mu_p\epsilon_{-1}+\mu_{p+1}\epsilon_1+\cdots+\mu_{p+q}\epsilon_q \in
P.$$ Then $W_{I(p,q)}$, which is isomorphic to the symmetric group
on $n$ letters, naturally acts on $\mu$. Given $\alpha=\sum_{i\in
[-p,q]}c_i\epsilon_i\in P$, put ${\bf
x}_{[-p,q]}^{\alpha}=\prod_{i\in[-p,q]}x_i^{c_i}$. Recall that the
Weyl formula for the Schur polynomial corresponding to $\mu$ is
given by
\begin{equation}\label{Weyl formula}
s_{\mu}({\bf x}_{[-p,q]})= \frac{\sum_{w\in
W_{I(p,q)}}(-1)^{\ell(w)}{\bf
x}_{[-p,q]}^{w(\mu+\delta_{p,q})-\delta_{p,q}}}{\prod_{i\in[-p]}\prod_{j\in
[q]}(1-x_i^{-1}x_j)}.
\end{equation}
Then we have a parabolic analogue as follows.

\begin{lem}\label{Weyl type} Suppose that $\mu\subset (d^n)$ for some $d>0$.
Following the above notations, we have
\begin{equation*}
{\bf x}_{[-p,q]}^{-d{\bf 1}^-_p}s_{\mu}({\bf x}_{[-p,q]})=
\frac{\sum_{w\in \W\cap W_{I(p,q)} } (-1)^{\ell(w)}
s_{(\lambda^{w,+})'}({\bf x}_{[q]})s_{(\lambda^{w,-})'}({\bf
x}_{[-p]}^{-1})}{\prod_{i\in[-p]}\prod_{j\in [q]}(1-x_i^{-1}x_j)},
\end{equation*}
where $\lambda=\mu'-(p^d)\in\mathbb{Z}_+^d$ and $\lambda^{w,\pm}$
are defined in \eqref{lambda w}.
\end{lem}
\pf By Lemma \ref{shifted action}, we have for $w\in\W\cap
W_{I(p,q)}$,
$$w(\mu+\delta_{p,q})-\delta_{p,q}-d{\bf 1}^-_p=-\sum_{i\in [-p]}\sigma_{-i}\epsilon_{i}+
\sum_{j\in [q]}\tau_{j}\epsilon_j,$$ for some $\sigma, \tau\in\cP$
with $\ell(\sigma)\leq p$ and $\ell(\tau)\leq q$. Now, for $w'\in
W_{I(p,q)\setminus\{0\}}$, we have
\begin{equation}
\begin{split}
&w'w(\mu+\delta_{p,q})-\delta_{p,q}-d{\bf 1}^-_p \\
&=w'\left(w(\mu+\delta_{p,q})-\delta_{p,q}-d{\bf 1}^-_p +\delta_{p,q}\right)-\delta_{p,q} \\
&=w'\left(-\sum_{i\in [-p]}\sigma_{-i}\epsilon_{i}+
\sum_{j\in [q]}\tau_{j}\epsilon_j+ q{\bf 1}^-_p+ \delta^-_{p}+ \delta^+_{q}\right)-q{\bf 1}^-_p- \delta^-_{p}- \delta^+_{q} \\
&= w'\left(-\sum_{i\in [-p]}\sigma_{-i}\epsilon_{i}+ \sum_{j\in
[q]}\tau_{j}\epsilon_j+  \delta^-_{p}+ \delta^+_{q}\right)-
\delta^-_{p}- \delta^+_{q},
\end{split}
\end{equation}
where $\delta^-_{p}=\sum_{i\in [-p]}(-i-1)\epsilon_{i}$,
$\delta^+_{q}=\sum_{j\in [q]}(q-j)\epsilon_j$ and
$\delta_{p,q}=q{\bf 1}^-_p+ \delta^-_{p}+ \delta^+_{q}$. Since
$s_{\nu}({\bf x}_{[-p]})=s_{\nu^*}({\bf x}_{[-p]}^{-1})$ for
$\nu\in\mathbb{Z}_+^p$, we obtain the result from \eqref{Weyl
formula}.\qed\vskip 3mm

In terms of $S_{\lambda}^{\A/\B}$, Lemma \ref{Weyl type} can be
written as follows.

\begin{lem}\label{Weyl type 1} For $\lambda\in\mathbb{Z}_+^d$, choose $p,q>0$ such that
$-p\leq \lambda_d\leq \lambda_1\leq q$. Then we have
$$S_{\lambda}^{[q]'/[-p]'}=
\frac{\sum_{w\in \W\cap W_{I(p,q)} } (-1)^{\ell(w)}
S_{\lambda^{w,+}}({\bf x}_{[q]'}) S_{\lambda^{w,-}}({\bf
x}_{[-p]'}^{-1})}{\prod_{i\in[-p]'}\prod_{j\in
[q]'}(1-x_i^{-1}x_j)}.
$$
\end{lem}
\pf Put $\mu=(\lambda+(p^d))'$. Note that $S_{\eta}({\bf
x}_{[q]'})=s_{\eta'}({\bf x}_{[q]'})$ and $S_{\eta}({\bf
x}^{-1}_{[-p]'})=s_{\eta'}({\bf x}^{-1}_{[-p]'})$ for $\eta\in\cP$.
Then, we have
\begin{equation}\label{rational form}
\begin{split}
{\bf x}_{[-p,q]'}^{-d{\bf 1}^-_p}s_{\mu}({\bf x}_{[-p,q]'})
&=\sum_{\nu}{\bf x}_{[-p,q]'}^{-d{\bf 1}^-_p}s_{\nu}({\bf x}_{[-p]'})s_{\mu/\nu}({\bf x}_{[q]'})  \\
&=\sum_{\nu}s_{\nu-(d^p)}({\bf
x}_{[-p]'})s_{\mu/\nu}({\bf x}_{[q]'})  \\
&=\sum_{\nu}s_{(d^p)/\nu}({\bf
x}^{-1}_{[-p]'})s_{\mu/\nu}({\bf x}_{[q]'})  \\
&=\sum_{\eta}S_{(p^d)/\eta}({\bf
x}^{-1}_{[-p]'})S_{\left(\lambda+(p^d)\right)/\eta}({\bf x}_{[q]'})  \\
&=S_{\lambda}({\bf x}_{[q]'};{\bf x}_{[-p]'}).
\end{split}
\end{equation}
The result follows from Lemma \ref{Weyl type} by replacing ${\bf
x}_{[-p,q]}$ with ${\bf x}_{[-p,q]'}$.

\qed

\begin{prop}\label{Weyl type 2}
For $\lambda\in\mathbb{Z}_+^d$, we have
\begin{equation*}
S_{\lambda}^{\mathbb{Z}_{>0}'/\mathbb{Z}_{< 0}'}= \frac{\sum_{w\in
\W} (-1)^{\ell(w)} S_{\lambda^{w,+}}({\bf x}_{\mathbb{Z}'_{>
0}})S_{\lambda^{w,-}}({\bf x}^{-1}_{\mathbb{Z}'_{< 0}})}
{\prod_{i,j}(1-x_i^{-1}x_j)},
\end{equation*}
and a Cauchy type identity
\begin{equation*}
\prod_{i,j,k}(1+x_i^{-1}z_k^{-1})(1+x_jz_k)
=\sum_{\lambda\in\mathbb{Z}_+^d}
S_{\lambda}^{\mathbb{Z}_{>0}'/\mathbb{Z}_{<0}'}s_{\lambda}({\bf
z}_{[d]}),
\end{equation*}
where $i\in\mathbb{Z}_{<0}'$, $j\in\mathbb{Z}_{>0}'$, $k\in [d]$,
and ${\bf z}_{[d]}=\{\,z_k\,|\,k\in [d]\,\}$.
\end{prop}
\pf First, it is easy to see that for $p,q>0$
$S_{\lambda}^{[q+1]'/[-p-1]'}$ reduces to $S_{\lambda}^{[q]'/[-p]'}$
when we put $x_{q+1}=0$ and $x^{-1}_{-p-1}=0$. Hence,
$S_{\lambda}^{[q]'/[-p]'}$ has the well-defined limit with respect
to both ${\bf x}_{[q]}$ and ${\bf x}^{-1}_{[-p]}$ when we let
$p,q\rightarrow \infty$, which is equal to
$S_{\lambda}^{\mathbb{Z}_{>0}'/\mathbb{Z}_{< 0}'}$. The first
identity follows from Lemma \ref{Weyl type 1}.

Next, consider the following dual Cauchy identity (cf.\cite{Mac}).
\begin{equation}
\prod_{i\in[n]}\prod_{j\in[d]}(1+x_iz_j)=\sum_{\mu\subset
(d^n)}s_{\mu}({\bf x}_{[n]})s_{\mu'}({\bf z}_{[d]}).
\end{equation}
Choose $p,q>0$ such that $p+q=n$. Replacing $[n]$ with $[-p,q]'$ and
multiplying $(x_{-p'}\cdots x_{-1'})^{-d}(z_1\cdots z_d)^{-p}$ on
both sides, we have
\begin{equation}\label{rationalizing}
\begin{split}
&\prod_{i\in[-p]'}\prod_{j\in[q]'}\prod_{k\in[d]}(1+x_i^{-1}z_k^{-1})(1+x_jz_k)\\
&=\sum_{\lambda\in\mathbb{Z}_+^d} {\bf x}^{-d{\bf
1}^{-1}_p}_{[-p]'}s_{\left(\lambda+(p^d)\right)'}({\bf
x}_{[-p,q]'})s_{\lambda}({\bf z}_{[d]}) \\
&=\sum_{\lambda\in\mathbb{Z}_+^d}
S_{\lambda}^{[q]'/[-p]'}s_{\lambda}({\bf z}_{[d]}). \ \ \ \
\text{(see \eqref{rational form})}
\end{split}
\end{equation}
 Hence, by letting $p,q\rightarrow \infty$,
we obtain the second identity. \qed

\begin{prop}\label{Weyl type 3} For $\lambda\in\mathbb{Z}_+^d$, we have
\begin{equation*}
S_{\lambda}^{\mathbb{Z}_{>0}/\mathbb{Z}_{< 0}}= \frac{\sum_{w\in \W}
(-1)^{\ell(w)} S_{\lambda^{w,+}}({\bf x}_{\mathbb{Z}_{>
0}})S_{\lambda^{w,-}}({\bf x}^{-1}_{\mathbb{Z}_{< 0}})}
{\prod_{i,j}(1-x_i^{-1}x_j)}, \\
\end{equation*}
and a Cauchy type identity
\begin{equation*}
\frac{1}{\prod_{i,j,k}(1-x_i^{-1}z_k^{-1})(1-x_jz_k)}
=\sum_{\lambda\in\mathbb{Z}_+^d}
S_{\lambda}^{\mathbb{Z}_{>0}/\mathbb{Z}_{< 0}}s_{\lambda}({\bf
z}_{[d]}),
\end{equation*}
where $i\in\mathbb{Z}_{<0}$, $j\in\mathbb{Z}_{>0}$, and $k\in [d]$.
\end{prop}
\pf Applying both $\omega_{\mathbb{Z}'_{< 0}}$ and
$\omega_{\mathbb{Z}'_{>0}}$ in Proposition \ref{Weyl type 2} and
then replacing ${\bf x}_{\mathbb{Z}'_{>0}}$ (resp. ${\bf
x}_{\mathbb{Z}'_{<0}}$) by ${\bf x}_{\mathbb{Z}_{>0}}$ (resp. ${\bf
x}_{\mathbb{Z}_{<0}}$), we obtain the identities. \qed

As a special case of Proposition \ref{Weyl type 3}, we obtain the
following Cauchy type identity with a restriction on the length of
partitions.
\begin{cor}\label{restricted Cauchy} For $d\geq 1$, we have
$$\sum_{\substack{\lambda\in \cP \\ \ell(\lambda)\leq d}}s_{\lambda}({\bf x})s_{\lambda}({\bf y})
=\frac{\sum_{w\in \W} (-1)^{\ell(w)} s_{{\bf 0}_d^{w,+}}({\bf
x})s_{{\bf 0}_d^{w,-}}({\bf y})}{\prod_{i,j\geq 1}(1-x_iy_j)},$$
where ${\bf x}={\bf x}_{\mathbb{Z}_{>0}}$, ${\bf y}={\bf
y}_{\mathbb{Z}_{>0}}$, and ${\bf
0}_d=(0,\ldots,0)\in\mathbb{Z}_+^d$.
\end{cor}
\pf Consider $S_{{\bf 0}_d}^{\mathbb{Z}_{>0}/\mathbb{Z}_{< 0}}$.
Replacing ${\bf x}^{-1}_{\mathbb{Z}_{< 0}}$ by ${\bf
y}_{\mathbb{Z}_{>0}}$ (that is, $x^{-1}_{-k}=y_{k}$ for $k\geq 1$),
we have
\begin{equation}
S_{{\bf 0}_d}^{\mathbb{Z}_{>0}/\mathbb{Z}_{< 0}}=
\sum_{\substack{\nu=(k^d)/\eta \\ k\geq 0,\ \eta\subset (k^d)
}}s_{\nu}({\bf x}) s_{\nu}({\bf y}) =\sum_{\ell(\lambda)\leq
d}s_{\lambda}({\bf x})s_{\lambda}({\bf y})
\end{equation}
since $s_{(k^d)/\eta}({\bf x})=s_{\lambda}({\bf x})$, where
$\lambda=(k^d)+\eta^*\in \cP$ (we regard $\eta$ as an element in
$\mathbb{Z}_+^d$) by the Littlewood-Richardson rule. Combining with
Proposition \ref{Weyl type 3}, we obtain the identity. \qed

\begin{rem}{\rm
Corollary  \ref{restricted Cauchy} can be stated more explicitly.
Given $w\in\W$, let $\mu=(\alpha|\beta)$ be the corresponding
partition (see Remark \ref{comments on cosets}) given in Frobenius
notation with $\delta(\mu)$ the length of $\alpha$ or $\beta$. Then
it is not difficult to see that
$${\bf 0}_d^{w,+}=\mu+(d^{\delta(\mu)}), \ \ {\bf 0}_d^{w,-}=\mu'+ (d^{\delta(\mu)}).$$
Hence, we obtain an alternative expression
\begin{equation}
\sum_{\substack{\lambda\in \cP \\ \ell(\lambda)\leq
d}}s_{\lambda}({\bf x})s_{\lambda}({\bf y}) =\frac{\sum_{\mu\in\cP}
(-1)^{|\mu|} s_{\mu+(d^{\delta(\mu)})}({\bf
x})s_{\mu'+(d^{\delta(\mu)})}({\bf y})}{\prod_{i,j\geq
1}(1-x_iy_j)}.
\end{equation}
 }
\end{rem}\vskip 3mm

\begin{thm}\label{main result} For $\lambda\in\mathbb{Z}_+^d$, we have
$$S_{\lambda}^{\A/\B}=
\frac{\sum_{w\in \W} (-1)^{\ell(w)} S_{\lambda^{w,+}}({\bf
x}_{\A})S_{\lambda^{w,-}}({\bf x}_{\B}^{-1})}{\Delta_{\A/\B}},$$
where
$$\Delta_{\A/\B}=\frac{\prod_{|a|=|b|}(1-x_ax_b^{-1})}{\prod_{|a|\neq |b|}(1+x_ax_b^{-1})}.$$
We also have the following Cauchy type identity
$$
\prod_{k\in [d]}\frac{\prod_{a\in \A_1}(1+x_az_k)\prod_{b\in
\B_1}(1+x_b^{-1}z_k^{-1})}{\prod_{a\in \A_0}(1-x_az_k)\prod_{b\in
\B_0}(1-x_b^{-1}z_k^{-1})}=
\sum_{\lambda\in\mathbb{Z}_+^d}S^{\A/\B}_{\lambda}s_{\lambda}({\bf
z}_{[d]}).
$$
\end{thm}
\pf For convenience, let us assume that $\A\subset \mathbb{Z}_{>0}$
and $\B\subset \mathbb{Z}_{<0}$ with arbitrary
$\mathbb{Z}_2$-gradings. Let $\A^{\circ}$ (resp. $\B^{\circ}$) be
the set of all positive (resp. negative) integers with a
$\mathbb{Z}_2$-grading such that $\A^{\circ}_i$ (resp.
$\B^{\circ}_i$) is infinite and $\A_i\subset\A^{\circ}_i$ (resp.
$\B_i\subset\B^{\circ}_i$) for $i\in\mathbb{Z}_2$.

We may view $\Lambda_{\mathbb{Z}_{>0}}$ (resp.
$\Lambda_{\mathbb{Z}_{< 0}}$) as a subring of
$\Lambda_{\A^{\circ}_0}\otimes \Lambda_{\A^{\circ}_1}$ (resp.
$\Lambda_{\B^{\circ}_0}\otimes \Lambda_{\B^{\circ}_1}$). Applying
$\omega_{\A^{\circ}_1}$ and $\omega_{\B^{\circ}_1}$ to
$S_{\lambda}^{\mathbb{Z}_{>0}/\mathbb{Z}_{< 0}}$ and the Cauchy type
identity in Proposition \ref{Weyl type 3} (we assume that the set of
variables in $\Lambda_{\mathbb{Z}_{< 0}}$ is ${\bf
x}^{-1}_{\mathbb{Z}_{< 0}}$), we obtain
\begin{equation}\label{Cauchy}
\begin{split}
&S_{\lambda}^{\A^{\circ}/\B^{\circ}}= \frac{\sum_{w\in \W}
(-1)^{\ell(w)} S_{\lambda^{w,+}}({\bf
x}_{\A^{\circ}})S_{\lambda^{w,-}}({\bf
x}_{\B^{\circ}}^{-1})}{\Delta_{\A^{\circ}/\B^{\circ}}}, \\
&\prod_{k\in [d]}\frac{\prod_{a\in
\A^{\circ}_1}(1+x_az_k)\prod_{b\in
\B_1^{\circ}}(1+x_b^{-1}z_k^{-1})}{\prod_{a\in
\A_0^{\circ}}(1-x_az_k)\prod_{b\in
\B_0^{\circ}}(1-x_b^{-1}z_k^{-1})}=
\sum_{\lambda\in\mathbb{Z}_+^d}S^{\A^{\circ}/\B^{\circ}}_{\lambda}s_{\lambda}({\bf
z}_{[n]}).
\end{split}
\end{equation}
Finally, by letting $x_a=x_b^{-1}=0$ for $a\in \A^{\circ}\setminus
\A$ and $b\in \B^{\circ}\setminus \B$ in \eqref{Cauchy}, we obtain
the results.\qed

\begin{rem}{\rm
The Cauchy type identity in Theorem \ref{main result} was also
proved in a bijective way in terms of tableaux \cite{Kwon}, which
can be viewed as an analogue of the Robinson-Schensted-Knuth
correspondence \cite{Kn}. }
\end{rem}

\subsection{Factorization property} Let us consider a particular case,
where we have a factorization property of $S^{\A/\B}_{\lambda}$.
Recall the following Cauchy type identity for skew Schur functions.

\begin{lem}[cf.\cite{Mac}] For $\lambda,\mu\in\cP$, we have
$$\sum_{\rho\in\cP}s_{\rho/\lambda}({\bf x})s_{\rho/\mu}({\bf y})
=\frac{1}{\prod_{i,j\geq
1}(1-x_iy_j)}\sum_{\tau\in\cP}s_{\mu/\tau}({\bf
x})s_{\lambda/\tau}({\bf y}),$$\qed
\end{lem}

By similar arguments as in Theorem \ref{main result}, it is
straightforward to rewrite the above identities as follows.
\begin{cor}\label{skewCauchy} For $\lambda,\mu\in\cP$, we have
$$\sum_{\rho\in\cP}S^{\A}_{\rho/\lambda}S^{\B}_{\rho/\mu}
=\frac{\prod_{|a|\neq|b|}(1+x_ax_b)}{\prod_{|a|=
|b|}(1-x_ax_b)}\sum_{\tau\in\cP}S^{\A}_{\mu/\tau}S^{\B}_{\lambda/\tau}.$$
\qed
\end{cor}

\begin{lem}\label{factorization of Schur}
Given $\lambda\in\mathbb{Z}_+^d$ and $m\in\mathbb{N}$, assume that
$d\geq m$ and $\lambda_m\geq 0$.  Let $p\geq 0$ and $\mu\in\cP$ be
such that $\lambda+(p^d)\in\cP$ and $\mu\subset \lambda+(p^d),
(p^d)$. Then
\begin{equation*}
s_{\left(\lambda+(p^d)\right)/\mu}({\bf x}_{[m]})=s_{\lambda^+}({\bf
x}_{[m]})s_{\nu/\mu}({\bf x}_{[m]}),
\end{equation*}
where $\nu=(p^d)+(\lambda^-)^*$.
\end{lem}
\pf Note that $s_{\rho/\tau}({\bf x}_{[m]})$ corresponding to a skew
Young diagram $\rho/\tau$ is the weight generating function of
$SST_{[m]}(\rho/\tau)$, the set of all semistandard tableaux of
shape $\rho/\tau$ with entries in $[m]$ (cf.\cite{Fu,Mac}).

Suppose that $SST_{[m]}\left(\left(\lambda+(p^d)\right)/\mu\right)$
is not empty. As usual, we enumerate the rows (resp. columns) in
$\lambda+(p^d)$ from top to bottom (resp. left to right). Also, we
may assume that $\nu/\mu$ is not empty, where
$\nu=(p^d)+(\lambda^-)^*$.

For $T\in SST_{[m]}\left(\left(\lambda+(p^d)\right)/\mu\right)$, let
$T_1$ (resp. $T_2$) be the subtableau obtained from the columns of
$T$ with indices greater than $p$ (resp. less than or equal to $p$).
The shapes of $T_1$ and $T_2$ are $\lambda^+$ and $\nu/\mu$,
respectively. So, this defines a map
\begin{equation}\label{factorization}
SST_{[m]}\left(\left(\lambda+(p^d)\right)/\mu\right)\longrightarrow
SST_{[m]}(\lambda^+) \times SST_{[m]}\left(\nu/\mu\right),
\end{equation}
by sending $T$ to $(T_1,T_2)$.

We claim that this is a one-to-one correspondence, which establishes
the corresponding identity of Schur polynomials. Let us construct an
inverse of the above map. Given a pair $(T_1,T_2)\in
SST_{[m]}(\lambda^+) \times SST_{[m]}\left(\nu/\mu\right)$, we
obtain a tableau $T$ (not necessarily semistandard) of shape
$\left(\lambda+(p^d)\right)/\mu$, where the first $p$ columns form
$T_2$, and the other columns form $T_1$. Note that the row and
column numbers are those in $\lambda+(p^d)$.

For $1\leq i\leq d$, let $a_i$ (resp. $b_i$) be the entry placed in
the $i$th row and the $p$th (resp. $(p+1)$th) column of $T$. We
assume that either $a_i$ or $b_i$ is empty if there is no entry.
Then $b_1,\ldots,b_d$ are the entries in the first column of $T_1$,
and $a_1,\ldots,a_d$ are the entries in the last column of $T_2$.

Since $m\leq d$, we have $b_{m+1}=\cdots=b_d=0$, and $b_k\geq k$ for
$1\leq k\leq m$. If all $a_i$ are empty, then it is clear that $T$
is semistandard. We assume that there exist non-empty entries
$a_s,\ldots, a_t$ with $1\leq s\leq t\leq d$. Since $\lambda_{m}\geq
0$, it follows that the $p$th column of $\lambda+(p^d)$ has at least
$m$ boxes, and hence $\left(\lambda+(p^d)\right)/\mu$ has at least
one box in rows lower than or equal to the $m$th row. This implies
that $t\geq m$, and $a_{k}\leq m+k-t\leq k\leq b_k$ for $s\leq k\leq
t$. Therefore, $T$ is a semistandard tableau of shape
$\left(\lambda+(p^d)\right)/\mu$, and the map $(T_1,T_2)\mapsto T$
is the inverse of \eqref{factorization}. This completes the proof.
\qed\vskip 3mm

\begin{thm}\label{factorization of SAB }
Suppose that $\A=[m]$ for $m>0$. For $\lambda\in\mathbb{Z}_+^d$, we
have
$$S^{[m]/\B}_{\lambda}=S_{\lambda^+}({\bf x}_{[m]})S_{\lambda^-}({\bf x}_{\B}^{-1})
 \Delta_{[m]/\B}^{-1}$$ if and only if $d\geq m$ and $\lambda_m\geq
0$.
\end{thm}
\pf First, note that by the Littlewood-Richardson rule, we have
$S^{\A}_{\mu}=S^{\A}_{\mu^{\pi}}$ for any $\A$ and $\mu\in\cP$,
where $\mu^{\pi}$ is the Young diagram obtained from $\mu$ by
$180^{\circ}$-rotation (cf.\cite{Mac}). Now, suppose that $d\geq m$
and $\lambda_m\geq 0$. Then we have  \allowdisplaybreaks{
\begin{equation}\label{factorization of SAB proof}
\begin{split}
&S^{[m]/\B}_{\lambda}
\\&=\sum_{\substack{\mu=\left(\lambda+(p^d)\right)/\eta,
 \\ \nu=(p^d)/\eta}}s_{\mu}({\bf x}_{[m]}) S_{\nu}({\bf
x}_{\B}^{-1}) \\
&=\sum_{\substack{\mu=\left(\lambda+(p^d)\right)/\eta,
 \\ \nu=(p^d)/\eta}}s_{\lambda^+}({\bf
x}_{[m]})s_{\left((p^d)+(\lambda^-)^*\right)/\eta}({\bf
x}_{[m]})S_{\nu}({\bf x}_{\B}^{-1}), \ \ \ \text{(by Lemma
\ref{factorization of Schur})} \\
&=s_{\lambda^+}({\bf x}_{[m]})\sum_{\substack{\nu=(p^d)/\eta
\\ \lambda^-\subset \nu^{\pi}}}s_{\nu^{\pi}/\lambda^-}({\bf
x}_{[m]})S_{\nu^{\pi}}({\bf x}_{\B}^{-1}) \ \ \ \ \ \ \ \ \ \ \ \ \ \ \ \ \ \text{(by $180^{\circ}$-rotation)} \\
&=s_{\lambda^+}({\bf x}_{[m]})\sum_{\lambda^-\subset
\tau}s_{\tau/\lambda^-}({\bf x}_{[m]})S_{\tau}({\bf x}_{\B}^{-1})
\\
&=S_{\lambda^+}({\bf x}_{[m]})S_{\lambda^-}({\bf x}_{\B}^{-1})
 \Delta_{[m]/\B}^{-1} \ \ \ \ \ \ \ \ \text{(by
Corollary \ref{skewCauchy})}.
\end{split}
\end{equation}}

Conversely, suppose that either $d< m$ or $\lambda_m< 0$. Let $p\geq
0$ and $\eta\in\cP$ be such that $\lambda+(p^d)\in\cP$ and
$\eta\subset \lambda+(p^d), (p^d)$. Then it is not difficult to see
that the difference
\begin{equation}
s_{\lambda^+}({\bf
x}_{[m]})s_{\left((p^d)+(\lambda^-)^*\right)/\eta}({\bf
x}_{[m]})-s_{\left(\lambda+(p^d)\right)/\eta}({\bf x}_{[m]}).
\end{equation}
is non-zero, in fact, an integral linear combination of monomials in
${\bf x}_{[m]}$ with non-negative coefficients. It implies that the
difference of $S^{[m]/\B}_{\lambda}$ and the last term in
\eqref{factorization of SAB proof} is non-zero. This completes the
proof.\qed

\begin{rem}{\rm
We can also prove Theorem \ref{factorization of SAB }   using
Theorem \ref{main result}. Suppose that $\lambda\in \mathbb{Z}_+^d$
is given. If $d\geq m$ and $\lambda_m\geq 0$, then it is not
difficult to see that $\ell(\lambda^{w,+})\leq m$ if and only if
$w=1 \in \W$, which implies that
$S^{[m]/\B}_{\lambda}=S_{\lambda^+}({\bf x}_{[m]})S_{\lambda^-}({\bf
x}_{\B}^{-1})\Delta_{[m]/\B}^{-1}$. Conversely, if $d< m$ or
$\lambda_m< 0$, then we can check that there exists at least one
non-trivial element $w\in \W$ such that $\ell(\lambda^{w,+})\leq m$,
which implies that
$S^{[m]/\B}_{\lambda}\Delta_{[m]/\B}-S_{\lambda^+}({\bf
x}_{[m]})S_{\lambda^-}({\bf x}_{\B}^{-1})$ is non-zero by Remark
\ref{comments on cosets} (2). }
\end{rem}

\subsection{Jacobi-Trudi formula} Finally, let us present another proof of the Jacobi-Trudi formula for
$S^{\A/\B}_{\lambda}$ \cite{Kwon} using the arguments given in 3.1.

Given $\A$ and $\B$, consider
\begin{equation}
H^{\A/\B}(t)=\frac{\prod_{a\in \A_1}(1+x_at)\prod_{b\in
\B_1}(1+x_b^{-1}t^{-1})}{\prod_{a\in \A_0}(1-x_at)\prod_{b\in
\B_0}(1-x_b^{-1}t^{-1})}=\sum_{k\in\mathbb{Z}}h^{\A/\B}_k t^k,
\end{equation}
where $h^{\A/\B}_k=S^{\A/\B}_k=\sum_{m-n=k}S_m({\bf x}_{\A})S_n({\bf
x}^{-1}_{\B})$ for $k\in\mathbb{Z}$.

\begin{prop}[cf.\cite{Kwon}]\label{JacobiTrudi}
For $\lambda\in\mathbb{Z}_+^d$, we have
$$S^{\A/\B}_{\lambda}={\rm det}(h^{\A/\B}_{\lambda_i-i+j})_{1\leq i,j\leq d}.$$
\end{prop}
\pf Choose $p,q>0$. Let $e_{k}({\bf x}_{[-p,q]'})$ be the $k$th
elementary symmetric polynomial in variables ${\bf x}_{[-p,q]'}$ for
$k\geq 0$, and $\widetilde{e}_k({\bf x}_{[-p,q]'})$ the coefficient
of $t^k$ for $-p\leq k \leq q$ in
$$\prod_{i\in[-p]'}\prod_{j\in[q]'}(1+x_i^{-1}t^{-1})(1+x_jt).$$
Then, from the classical Jacobi-Trudi formula, we have
\begin{equation}
\begin{split}
S_{\lambda}^{[q]'/[-p]'}&=  (x_{-p'}\cdots
x_{-1'})^{-d}s_{\left(\lambda+(p^d)\right)'}({\bf x}_{[-p,q]'})
\\
&= (x_{-p'}\cdots x_{-1'})^{-d}{\rm
det}(e_{\lambda_i+p-i+j}({\bf x}_{[-p,q]'}))_{1\leq i,j\leq d} \\
&={\rm det}(\widetilde{e}_{\lambda_i-i+j}({\bf x}_{[-p,q]'}))_{1\leq
i,j\leq d}.
\end{split}
\end{equation}
If we follow the same arguments as in Theorem \ref{main result},
then $S_{\lambda}^{[q]'/[-p]'}$ and $\widetilde{e}_k({\bf
x}_{[-p,q]'})$ are replaced by $S_{\lambda}^{\A/\B}$ and
$h^{\A/\B}_k$, respectively. This completes the proof. \qed

\section{Applications}

In this section, we discuss applications of our Weyl type formula
for $S^{\A/\B}_{\lambda}$ in representation theory. We assume that
the ground field is $\mathbb{C}$.

\subsection{Representations of $\frak{gl}_{m|n}$}
First, let us derive a Weyl type formula for hook Schur polynomials.
Though its proof is purely combinatorial, we will translate the
result in the language of representation theory to see its relation
with   finite dimensional irreducible representations of
$\frak{gl}_{m|n}$. So, let us give  brief review on representations
of the Lie superalgebra $\frak{gl}_{m|n}$ (cf.\cite{Kac77}).

For non-negative integers $m$ and $n$, not both zero, let
$\mathbb{C}^{m|n}=\mathbb{C}^{m|0}\oplus\mathbb{C}^{0|n}$ be the
$(m+n)$-dimensional superspace with the even subspace
$\mathbb{C}^{m|0}=\mathbb{C}^m$ and the odd subspace
$\mathbb{C}^{0|n}=\mathbb{C}^n$. We denote by
$\{\,\epsilon_i\,|\,i\in [-m]\,\}$ and $\{\,\epsilon_j\,|\,j\in
[n]'\,\}$ the homogeneous bases of $\mathbb{C}^{m|0}$ and
$\mathbb{C}^{0|n}$ respectively, which form a standard basis of
$\mathbb{C}^{m|n}$. Then the space of $\mathbb{C}$-linear
endomorphisms of $\mathbb{C}^{m|n}$ is naturally equipped with a
$\mathbb{Z}_2$-grading, and becomes a Lie superalgebra with respect
to a super bracket, which is  called a general linear superalgebra
$\frak{gl}_{m|n}$. Put $[-m,n']=[-m]\cup [n]'$.

We may identify $\frak{g}=\frak{gl}_{m|n}$ with the set of
$(m+n)\times (m+n)$ matrices with respect to the standard basis of
$\mathbb{C}^{m|n}$. Then the subspace $\mathfrak{h}$ of diagonal
matrices forms a Cartan subalgebra, and under the adjoint action of
$\mathfrak{h}$ on $\mathfrak{g}$, we have a root space
decomposition, $\frak{g}=\frak{h}\oplus \left(\bigoplus_{\alpha\in
\Delta} \frak{g}_{\alpha}\right)$, where $\Delta$ is the set of all
roots of $\frak{g}$. Let $\Delta^{+}$ be the set of positive roots,
and $\Delta^+_0$ (resp. $\Delta^+_1$) the set of positive even
(resp. odd) roots. Since $\frak{h}^*=\mathbb{C}^{m|n}$, we have
\begin{equation}
\begin{split}
\Delta^+_0&=\{\,\epsilon_{i}-\epsilon_j\,|\,i,j\in[-m,n'],\
|i|=|j|,\
i<j\,\}, \\
\Delta^+_1&=\{\,\epsilon_{i}-\epsilon_j\,|\,i\in [-m],\ j\in
[n]'\,\},
\end{split}
\end{equation}
where we assume that $-m<-m+1<\cdots<-1<1'<2'\cdots<n'$. We consider
a $\mathbb{Z}$-grading $\frak{g}=\frak{g}_{-1}\oplus\frak{g}_0\oplus
\frak{g}_1$ consistent with its parity, where
$\frak{g}_0=\frak{gl}_m\oplus \frak{gl}_n$, and a triangular
decomposition, $\mathfrak{g}= \frak{n^-}\oplus\mathfrak{h}\oplus
\frak{n^+}$ such that $ \frak{g}_{\pm 1}\subset\frak{n}^{\pm}$.

Let ${P}_{m|n}$ be the $\mathbb{Z}$-lattice of $\frak{h}^*$
generated by $\{\,\epsilon_i\,|\,i\in [-m,n']\,\}$, which is called
the set of integral weights, and let ${P}_{m|n}^+$ be the set of
weights $\Lambda=\sum_{i\in[-m,n']}\Lambda_i\epsilon_i\in {P}_{m|n}$
such that $\Lambda_{-m}\geq \ldots\geq \Lambda_{-1}$ and
$\Lambda_{1'}\geq \ldots\geq \Lambda_{n'}$. An element in
${P}_{m|n}^+$ is called a dominant integral weight.

Suppose that $M$ is a finite-dimensional $\frak{gl}_{m|n}$-module,
which is $\frak{h}$-diagonalizable. Then, we have a  weight
decomposition $M=\bigoplus_{\lambda\in \frak{h}^*}M_{\lambda}$. When
$\dim M_{\lambda}\neq 0$, we call $\lambda$ a  weight of $M$. For
convenience, we assume that all the weight are integral. We define
the character of $M$ by ${\rm ch}M=\sum_{\lambda\in P_{m|n}}{\rm
dim}M_{\lambda}e^{\lambda}$, where $\{\,e^{\lambda}\,|\,\lambda\in
{P}_{m|n}\,\}$ is the set of formal variables.

Given $\Lambda\in {P}_{m|n}^+$, let $L^0(\Lambda)$ be the finite
dimensional irreducible highest weight $\frak{g}_0$-module with
highest weight $\Lambda$. We may view $L^0(\Lambda)$ as a
$\frak{g}_0\oplus \frak{g}_1$-module, where $\frak{g}_1$ acts
trivially on $L^0(\Lambda)$. The Kac module $K_{m|n}(\Lambda)$ is
defined to be the induced representation
$K_{m|n}(\Lambda)=U(\frak{g})\otimes_{U(\frak{g}_0\oplus
\frak{g}_1)}L^0(\Lambda)$ where $U(\frak{g})$ and
$U(\frak{g}_0\oplus \frak{g}_1)$ are the enveloping algebras, and it
has a unique maximal irreducible quotient $L_{m|n}(\Lambda)$. Then
$\{\,L_{m|n}(\Lambda)\,|\,\Lambda\in {P}_{m|n}^+\,\}$ forms a
complete set of pairwise non-isomorphic finite dimensional
irreducible representations of $\frak{gl}_{m|n}$ with integral
weights.

In \cite{Kac77-2,Kac78}, Kac gave a necessary and sufficient
condition for $K_{m|n}(\Lambda)$ ($\Lambda\in {P}_{m|n}^+$) to be
irreducible, and called such weights typical. For an  atypical
weight $\Lambda\in {P}_{m|n}^+$, which is not typical,
$L_{m|n}(\Lambda)$ has a resolution, where each term has a
filtration with quotients isomorphic to Kac modules (see \cite{S}),
and its character is given by
\begin{equation}\label{character formula}
{\rm ch}L_{m|n}(\Lambda)=\sum_{\Lambda'\in
P^+_{m|n}}a_{\Lambda\,\Lambda'} {\rm ch}K_{m|n}(\Lambda'),
\end{equation}
for some $a_{\Lambda\,\Lambda'}\in\mathbb{Z}$. In \cite{S},
Serganova gave an algorithm for computing these coefficients using
the geometry of the associated supergroups. Recently, in \cite{Br},
Brundan gave another algorithm using a remarkable connection with
canonical bases of quantum group $U_q(\frak{gl}_{\infty})$.\vskip
3mm

Now, let us give an explicit expression of \eqref{character formula}
for irreducible tensor representations. For
$\Lambda=\sum_{i\in[-m,n']}\Lambda_i\epsilon_i\in {P}_{m|n}^+$, we
may identify $\Lambda$ with a pair of generalized partitions given
by
\begin{equation}\label{bipartition}
\Lambda^{<0}=(\Lambda_{-m},\ldots,\Lambda_{-1})\in\mathbb{Z}_+^m, \
\ \
\Lambda^{>0}=(\Lambda_{1'},\ldots,\Lambda_{n'})\in\mathbb{Z}_+^n.
\end{equation}
If we put $x_i=e^{\epsilon_i}$ for $i\in [-m]$, and
$y_j=e^{\epsilon_{j'}}$ for $j\in [n]$, then the character of
$K_{m|n}(\Lambda)$ is given by
\begin{equation}
{\rm ch}K_{m|n}(\Lambda)=s_{\Lambda^{<0}}({\bf
x}_{[-m]})s_{\Lambda^{>0}}({\bf
y}_{[n]})\prod_{\substack{i\in [-m], \\
j\in [n]}}(1+x_{i}^{-1}y_j).
\end{equation}
Let $\mathcal{P}_{m|n}^{+}$ be the set of dominant integral weights
$\Lambda$ such that $\Lambda_k\geq 0$ for all $k\in [-m,n']$ and
$\ell((\Lambda^{>0})')\leq \Lambda_{-1}$. For
$\Lambda\in\mathcal{P}_{m|n}^+$, we define a partition
\begin{equation}\label{hook partition}
\lambda(\Lambda)=(\Lambda_{-m},\ldots,\Lambda_{-1},\Lambda'_{1},\ldots,\Lambda'_{\ell}),
\end{equation}
where $(\Lambda^{>0})'=(\Lambda'_{1},\ldots,\Lambda'_{\ell})$ and
$\ell=\ell((\Lambda^{>0})')$. Then the map $\Lambda \mapsto
\lambda(\Lambda)$ gives a one-to-one correspondence between
$\mathcal{P}_{m|n}^+$ and
$\cP_{m|n}=\{\,\lambda\in\cP\,|\,\lambda_{m+1}\leq n\,\}$, the set
of all $(m,n)$-hook partitions.

In \cite{BR}, it was shown that the tensor algebra
$\mathscr{T}(\mathbb{C}^{m|n})$ generated by the natural
representation $\mathbb{C}^{m|n}$ is completely reducible, and for
$\Lambda\in {P}^+_{m|n}$, $L_{m|n}(\Lambda)$  occurs in
$\mathscr{T}(\mathbb{C}^{m|n})$ if and only if $\Lambda\in
\mathcal{P}_{m|n}^+$. Moreover, we have
\begin{equation}\label{hook Schur polynomial}
{\rm ch}L_{m|n}(\Lambda)=\sum_{\mu\subset \lambda}s_{\mu}({\bf
x}_{[-m]})s_{\lambda'/\mu'}({\bf y}_{[n]})
\end{equation}
with $\lambda=\lambda(\Lambda)\in \cP_{m|n}$, which is called the
hook Schur polynomial corresponding to $\lambda$. Following our
notation, ${\rm
ch}L_{m|n}(\Lambda)=S_{\lambda}^{[-m,n']}=S_{\lambda}({\bf
x}_{[-m,n']})$, where we identify ${\bf x}_{[n]'}\subset {\bf
x}_{[-m,n']}$ with ${\bf y}_{[n]}$.

\begin{lem}\label{Weyl type 4}
For $\Lambda\in \mathcal{P}_{m|n}^+$, let
$\lambda=\lambda(\Lambda)\in\cP_{m|n}$ and $\nu=\lambda'-(m^d)$ for
$d\geq \lambda_1$. Then $\nu\in\mathbb{Z}_+^d$, and
\begin{equation*}
(x_{-m}\cdots x_{-1})^{-d}{\rm ch}L_{m|n}(\Lambda) = S_{\nu}({\bf
x}_{[n]};{\bf x}_{[-m]'}),
\end{equation*}
where we identify   ${\bf x}_{[-m]}$ with ${\bf x}_{[-m]'}$, and
${\bf y}_{[n]}$ with ${\bf x}_{[n]}$.
\end{lem}
\pf By similar arguments as in Lemma \ref{Weyl type 1}
\eqref{rational form}, we have
{\allowdisplaybreaks
\begin{equation*}
\begin{split}
&(x_{-m}\cdots x_{-1})^{-d} {\rm ch}L_{m|n}(\Lambda) \\
&=(x_{-m}\cdots x_{-1})^{-d}\left( \sum_{\mu\subset
\lambda}s_{\mu}({\bf
x}_{[-m]})s_{\lambda'/\mu'}({\bf y}_{[n]}) \right) \\
&=\sum_{\mu\subset \lambda}s_{\mu-(d^m)}({\bf
x}_{[-m]})s_{\lambda'/\mu'}({\bf y}_{[n]}) \\
&=\sum_{\mu\subset \lambda}s_{(d^m)/\mu}({\bf
x}^{-1}_{[-m]})s_{\lambda'/\mu'}({\bf y}_{[n]})  \\
&=\sum_{\eta}S_{(m^d)/\eta}({\bf
x}^{-1}_{[-m]'})S_{\left(\nu+(m^d)\right)/\eta}({\bf y}_{[n]})
\ \ \ \ \text{by replacing ${\bf x}_{[-m]}$ with ${\bf x}_{[-m]'}$} \\
&=S_{\nu}({\bf x}_{[n]};{\bf x}_{[-m]'})=S_{\nu}^{[n]/[-m]'} \ \ \ \
\ \ \ \ \ \ \ \ \ \ \text{by replacing ${\bf y}_{[n]}$ with ${\bf
x}_{[n]}$}.
\end{split}
\end{equation*}}
\qed

We have another proof of the factorization property of hook Schur
polynomials, or the irreducibility of the Kac-module
$K_{m|n}(\Lambda)$ for $\Lambda\in \mathcal{P}^+_{m|n}$.

\begin{cor}[\cite{BR,rem}] For $\Lambda\in \mathcal{P}^+_{m|n}$,
we have $L_{m|n}(\Lambda)=K_{m|n}(\Lambda)$ if and only if
$\Lambda_{-1}\geq n$.
\end{cor}
\pf For $\Lambda\in \mathcal{P}^+_{m|n}$, $\Lambda_{-1}\geq n$ if
and only if $\lambda=\lambda(\Lambda)\in \cP_{m|n}$ contains the
partition $(n^m)$ of a rectangular shape. By Lemma \ref{Weyl type
4}, ${\rm ch}L_{m|n}(\Lambda)=(x_{-m}\cdots x_{-1})^{d}
S^{[n]/[-m]'}_{\nu}$, where $\nu=\lambda'-(m^d)$ for some $d\geq
\lambda_1$. Hence $\Lambda_{-1}\geq n$ if and only if $\nu$
satisfies the condition in Theorem \ref{factorization of SAB }, when
we have ${\rm ch}L_{m|n}(\Lambda)={\rm ch}K_{m|n}(\Lambda)$. This
completes the proof. \qed

For $\Lambda\in \mathcal{P}_{m|n}^+$ with
$\lambda=\lambda(\Lambda)$, set
\begin{equation}
\mathcal{W}_{m|n}(\Lambda)=\{\,w\in\mathcal{W}\,|\,\ell((\nu^{w,-})')\leq
m,\ \ell(\nu^{w,+})\leq n\,\},
\end{equation}
where $\nu=\lambda'-(m^d)\in\mathbb{Z}_+^d$ with $d\geq \lambda_1$.
For $w\in \mathcal{W}_{m|n}(\Lambda)$, we define
\begin{equation}
w\ast\Lambda= (\left((\nu^{w,-})'-(d^m)\right)^*,\nu^{w,+})\in
{P}_{m|n}^+,
\end{equation}
where we identify ${P}_{m|n}^+$ with $\mathbb{Z}_+^m\times
\mathbb{Z}_+^n$ (cf.\eqref{bipartition}). We can check that given
$w\in \mathcal{W}_{m|n}(\Lambda)$, $w\in W_{I(m,q)}$ for some $q$,
and $w\ast \Lambda$ does not depend on the choice of $d$. Now, we
obtain a Weyl type formula for hook Schur polynomials, which
recovers the Cheng and Zhang's formula \cite{CZ} in a new
combinatorial way.

\begin{thm}[\cite{CZ}]\label{Weyl for hook Schur} For $\Lambda\in \mathcal{P}^+_{m|n}$, we have
\begin{equation*}
{\rm ch}L_{m|n}(\Lambda)=\sum_{w\in
\mathcal{W}_{m|n}(\Lambda)}(-1)^{\ell(w)}{\rm
ch}K_{m|n}(w\ast\Lambda).
\end{equation*}
\end{thm}
\pf By Lemma \ref{Weyl type 4} and Theorem \ref{main result}, we
have {\allowdisplaybreaks
\begin{equation*}
\begin{split}
&{\rm ch}L_{m|n}(\Lambda) \\ &=(x_{-m}\cdots x_{-1})^d
\left(\sum_{w\in \W} (-1)^{\ell(w)} s_{\nu^{w,+}}({\bf
y}_{[n]})s_{(\nu^{w,-})'}({\bf
x}_{[-m]}^{-1})\right)\prod_{\substack{i\in [-m] \\
j\in [n]}}(1+x_i^{-1}y_j) \\
&=\left(\sum_{w\in \W} (-1)^{\ell(w)} s_{\nu^{w,+}}({\bf
y}_{[n]})s_{(\nu^{w,-})'-(d^m)}({\bf
x}_{[-m]}^{-1})\right)\prod_{\substack{i\in [-m] \\
j\in [n]}}(1+x_i^{-1}y_j) \\
&=\left(\sum_{w\in \W} (-1)^{\ell(w)} s_{\nu^{w,+}}({\bf
y}_{[n]})s_{\left((\nu^{w,-})'-(d^m)\right)^*}({\bf
x}_{[-m]})\right)\prod_{\substack{i\in [-m] \\
j\in [n]}}(1+x_i^{-1}y_j) \\
&=\sum_{w\in \mathcal{W}_{m|n}(\Lambda)}(-1)^{\ell(w)}{\rm
ch}K(w\ast\Lambda).
\end{split}
\end{equation*}}
\qed

\begin{rem}{\rm
Lemma \ref{Weyl type 4} and hence Theorem \ref{Weyl for hook Schur}
can also be proved from a Howe duality of
$(\frak{gl}_{m|n},\frak{gl}_d)$ acting on the supersymmetric algebra
generated by $(\mathbb{C}^{m|0*}\oplus \mathbb{C}^{0|n})\oplus
\mathbb{C}^d$ (cf.\cite{CLZ,CW01}) without a combinatorial
definition of ${\rm ch}L_{m|n}(\Lambda)$ by Berele and Regev.}
\end{rem}

The denominator identity for $\frak{gl}_{m|n}$ is given as follows.
\begin{cor} For $w\in \mathcal{W}_{m|n}(0)$, let $w\ast
0=(\theta^{w}_m,\theta^{w}_n)\in {P}_{m|n}^+$. Then
\begin{equation*}
\begin{split}
&\prod_{i,j}(x_i+y_j)^{-1}\prod_{i<i'}(x_i-x_{i'})
\prod_{j<j'}(y_j-y_{j'}) \\
&=\sum_{w\in \mathcal{W}_{m|n}(0)} \sum_{\substack{w_1\in S_m
\\ w_2\in S_n}}(-1)^{\ell(w_1)+\ell(w_2)+\ell(w)}
{\bf x}_{[-m]}^{w_1(\theta^{w}_m+\delta_m)-n{\bf 1}_m} {\bf
y}_{[n]}^{w_2(\theta^{w}_n+\delta_n)},
\end{split}
\end{equation*}
where $i,i'\in [-m]$, $j,j'\in [n]$, $\delta_m=\sum_{i\in
[-m]}(-i-1)\epsilon_i$,$\delta_n=\sum_{j\in [n]}(j-1)\epsilon_j$,
${\bf 1}_m=\sum_{i\in [-m]}\epsilon_i$, and $S_m$ {\rm(}resp.
$S_n${\rm )} is the symmetric group on the letters $[-m]$
{\rm(}resp. $[n]${\rm)}.
\end{cor}
\pf Since the hook Schur polynomial corresponding to empty partition
(that is, $\Lambda=0$) is $1$, the identity follows from Theorem
\ref{Weyl for hook Schur}. \qed \vskip 3mm

\begin{rem}{\rm
For $w\in \mathcal{W}_{m|n}(0)$, let $\mu$ be the corresponding
partition (see Remark \ref{comments on cosets}). Then  $w\ast 0 \in
{P}_{m|n}^+$ if and only if $\ell(\mu)\leq m,n$. Moreover, in this
case, we have
\begin{equation*}
w\ast
0=(\theta^{w}_m,\theta^{w}_n)=(-\mu_m,\ldots,-\mu_1,\mu_1,\ldots,\mu_n)\in
{P}_{m|n}^+.
\end{equation*}
Since $|\mu|=\ell(\mu)$, we may write by abuse of notation
\begin{equation*}
\begin{split}
&\prod_{i,j}(x_i+y_j)^{-1}\prod_{i<i'}(x_i-x_{i'})
\prod_{j<j'}(y_j-y_{j'}) \\
&=\sum_{\ell(\mu)\leq m,n}(-1)^{|\mu|} \sum_{\substack{w_1\in S_m
\\ w_2\in S_n}}(-1)^{\ell(w_1)+\ell(w_2)}
{\bf x}_{[-m]}^{w_1(\mu^*+\delta_m)-n{\bf 1}_m} {\bf
y}_{[n]}^{w_2(\mu+\delta_n)}.
\end{split}
\end{equation*}
 }
\end{rem}

\subsection{Representations of $\frak{gl}_{m+n}$}
Let us consider a character formula  for a certain class of infinite
dimensional highest weight representations of the Lie algebra
$\frak{gl}_{m+n}$, which comes from the study of unitary highest
weight representations of its associated Lie group $U(m,n)$
(cf.\cite{EHW,H,KV}) or from a parabolic analogue of Kazhdan-Lusztig
theory (cf.\cite{CC,Deo}).

For non-negative integers $m$ and $n$, not both zero, let
$\mathbb{C}^{m+n}=\mathbb{C}^{m}\oplus\mathbb{C}^{n}$ be the
$(m+n)$-dimensional space, where $\{\,\epsilon_i\,|\,i\in [-m]\,\}$
(resp. $\{\,\epsilon_j\,|\,j\in [n]\,\}$) is the basis of
$\mathbb{C}^{m}$ (resp. $\mathbb{C}^{n}$). We may identify the
general linear algebra $\frak{g}=\frak{gl}_{m+n}$ with the set of
$(m+n)\times (m+n)$ matrices whose row and column indices are from
$[-m]\cup [n]$ (or simply $[-m,n]$).

We denote by $\frak{h}$ and $\frak{b}$ the Cartan subalgebra of the
diagonal matrices and the Borel subalgebra of the upper triangular
matrices respectively. We also put $\frak{t}=\frak{gl}_m\oplus
\frak{gl}_n$, which is naturally embedded in $\frak{g}$. Let
$P_{m+n}$ be the $\mathbb{Z}$-lattice of $\frak{h}^*$ generated by
$\{\,\epsilon_i\,|\,i\in [-m,n]\,\}$, and   ${P}_{m+n}^+$   the set
of $\frak{t}$-dominant integral weights, that is, the set of weights
$\Lambda=\sum_{i\in[-m,n]}\Lambda_i\epsilon_i\in {P}_{m+n}$ such
that $\Lambda_{-m}\geq \ldots\geq \Lambda_{-1}$ and $\Lambda_{1}\geq
\ldots\geq \Lambda_{n}$.

Given $\Lambda\in {P}_{m+n}^+$, let $L^0(\Lambda)$ be the finite
dimensional irreducible highest weight $\frak{t}$-module with
highest weight $\Lambda$. We may view $L^0(\Lambda)$ as a
representation of the parabolic subalgebra $\frak{q}=\frak{t} +
\frak{b}$, where the action is extended in a trivial way. Now we
define the  generalized Verma module $V_{m+n}(\Lambda)$ to be the
induced representation
$V_{m+n}(\Lambda)=U(\frak{g})\otimes_{U(\frak{q})}L^0(\Lambda)$.
Then $V_{m+n}(\Lambda)$ has a unique maximal irreducible quotient
$L_{m+n}(\Lambda)$. Similarly, we define the characters of
$L_{m+n}(\Lambda)$ and $V_{m+n}(\Lambda)$ in terms of the formal
variables $e^{\lambda}$ ($\lambda\in P_{m+n}$). Put
$x_i=e^{\epsilon_i}$ for $i\in [-m,n]$. Then
$${\rm ch}V_{m+n}(\Lambda)=\frac{s_{\Lambda^{<0}}({\bf
x}_{[-m]})s_{\Lambda^{>0}}({\bf x}_{[n]})}{\prod_{i,j}
(1-x_i^{-1}x_j)},$$ where $i\in [-m]$, $j\in [n]$,
$\Lambda^{<0}=(\Lambda_{-m},\ldots,\Lambda_{-1})$ and
$\Lambda^{>0}=(\Lambda_{1},\ldots,\Lambda_{n})$.

From the Kazhdan-Lusztig conjecture \cite{KL} proved in \cite{BB,BK}
and its parabolic analogue \cite{CC,Deo}, we have
\begin{equation}
{\rm ch}L_{m+n}(\Lambda)=\sum_{\Lambda'\in
P^+_{m+n}}b_{\Lambda\,\Lambda'}{\rm ch}V_{m+n}(\Lambda'),
\end{equation}
where the integers $b_{\Lambda\,\Lambda'}$ are determined explicitly
in terms of (parabolic) Kazhdan-Lusztig polynomials evaluated at
$1$.

Now, using Howe duality let us derive an explicit expression of
${\rm ch}L_{m+n}(\Lambda)$ for particular highest weights in
$P^+_{m+n}$, which is equivalent to the Enright's formula given in
\cite{En}. For $d\in\mathbb{N}$, let $P_d=\bigoplus_{k\in
[d]}\mathbb{Z}\varepsilon_k$ be the weight lattice of $\frak{gl}_d$.
Consider the symmetric algebra $\mathscr{S}_d$ generated by
$(\mathbb{C}^{m^*}\otimes \mathbb{C}^{d^*})\oplus
(\mathbb{C}^{n}\otimes \mathbb{C}^{d})$, where $\mathbb{C}^{m^*}$
and $\mathbb{C}^{d^*}$ are the duals of the natural representations
of $\frak{gl}_m$ and $\frak{gl}_d$ respectively. Then there is a
semi-simple $(\frak{gl}_{m+n},\frak{gl}_d)$-action on
$\mathscr{S}_d$, which gives the following multiplicity free
decomposition \cite{H,KV}
\begin{equation}\label{Howeduality}
\mathscr{S}_d=\bigoplus_{\lambda\in\mathbb{Z}_+^d}L_{m+n}(\Lambda(\lambda))\otimes
L_d(\lambda),
\end{equation}
where $L_d(\lambda)$ is the irreducible highest weight
$\frak{gl}_d$-module with highest weight
$\lambda=(\lambda_1,\ldots,\lambda_d)\in\mathbb{Z}_+^d$ (or
$\lambda=\sum_{k\in [d]}\lambda_k\varepsilon_k\in P_d$), and
$\Lambda(\lambda)$ is a highest weight in $P_{m+n}^+$. In fact, if
\begin{equation}
\lambda=\lambda^+ +
(\lambda^-)^*=(\lambda^+_1,\ldots,\lambda^+_p,0,\ldots,0,-\lambda^-_q,\ldots,-\lambda^-_1),
\end{equation}
then
\begin{equation}
\Lambda(\lambda)=(\underbrace{-d,\ldots,-d,-\lambda^-_q-d,\ldots,-\lambda^-_1-d}_{m},
\underbrace{\lambda^+_1,\ldots,\lambda^+_p,0,\ldots,0}_{n}),
\end{equation}
where we identify $P_{m+n}$ with $\mathbb{Z}^{m+n}$.

\begin{prop}\label{Weyl type 5} For $\lambda\in\mathbb{Z}_+^d$, we have
$$(x_{-m}\cdots
x_{-1})^{d}{\rm ch}L_{m+n}(\Lambda(\lambda))=S_{\lambda}({\bf
x}_{[n]};{\bf x}_{[-m]}).$$
\end{prop}
\pf  The decomposition in \eqref{Howeduality} gives the following
identity
\begin{equation}
\frac{(x_{-m}\cdots
x_{-1})^{-d}}{\prod_{i,j,k}(1-x_iz_k)(1-x_j^{-1}z_k^{-1})}=
\sum_{\lambda\in\mathbb{Z}_+^d}{\rm
ch}L_{m+n}(\Lambda(\lambda))s_{\lambda}({\bf z}_{[d]}),
\end{equation}
where $i\in [n]$, $j\in [-m]$, $k\in [d]$, and ${\bf
z}_{[d]}=\{z_k=e^{\varepsilon_k}\,|\,k\in [d]\,\}$. Comparing with
the Cauchy type identity in Theorem \ref{main result}, it follows
from the linear independence of rational Schur polynomials that
$(x_{-m}\cdots x_{-1})^{d}{\rm
ch}L_{m+n}(\Lambda(\lambda))=S_{\lambda}^{[n]/[-m]}$. \qed

\begin{cor} For  $\lambda\in\mathbb{Z}_+^d$, we have
$L_{m+n}(\Lambda(\lambda))=V_{m+n}(\Lambda(\lambda))$ if and only if
$d\geq n$ and $\lambda_n\geq 0$.
\end{cor}
\pf It  follows from Theorem \ref{factorization of SAB }. \qed\vskip
3mm

We also have an interesting analogue of the Jacobi-Trudi formula
from Proposition \ref{JacobiTrudi}.
\begin{cor} For  $\lambda\in\mathbb{Z}_+^d$, we have
$${\rm ch}L_{m+n}(\Lambda(\lambda))={\rm det}({\rm ch}L_{m+n}(\Lambda(\lambda_i-i+j)))_{1\leq i,j\leq
d}.$$\qed
\end{cor}

Suppose that $\Lambda\in {P}_{m+n}^+$ is given and
$\Lambda=\Lambda(\lambda)$ for some $\lambda\in\mathbb{Z}_+^d$. Set
\begin{equation}
\mathcal{W}_{m+n}(\Lambda)=\{\,w\in\mathcal{W}\,|\,\ell(\lambda^{w,-})\leq
m,\ \ell(\lambda^{w,+})\leq n\,\}.
\end{equation}
For $w\in \mathcal{W}_{m+n}(\Lambda)$, we define
\begin{equation}
w\ast\Lambda= \left( (\lambda^{w,-}+(d^m))^* ,
\lambda^{w,+}\right)\in {P}_{m+n}^+,
\end{equation}
where we identify  ${P}_{m+n}^+$ with $\mathbb{Z}_+^m\times
\mathbb{Z}_+^n$. Now, we can state a Weyl type formula for ${\rm
ch}L_{m+n}(\Lambda(\lambda))$ ($\lambda\in\mathbb{Z}_+^d$).
\begin{thm}[cf.\cite{En}]\label{Weyl for unitary}
Given $\Lambda\in {P}^+_{m+n}$ with $\Lambda=\Lambda(\lambda)$ for
some $\lambda\in\mathbb{Z}_+^d$, we have
\begin{equation*}
{\rm ch}L_{m+n}(\Lambda)=\sum_{w\in
\mathcal{W}_{m+n}(\Lambda)}(-1)^{\ell(w)}{\rm
ch}V_{m+n}(w\ast\Lambda).
\end{equation*}
\end{thm}
\pf It follows from Proposition \ref{Weyl type 5} with Theorem
\ref{main result}. The proof is similar to that of Theorem \ref{Weyl
for hook Schur}. \qed

\begin{rem}{\rm
The parametrization of highest weights for generalized Verma modules
in Theorem \ref{Weyl for unitary} is different from the one in
\cite{En}, where  the sum is given over the Weyl group of
$\frak{gl}_{m+n}$ with its shifted action on the highest weight
$\Lambda$, which is not equal to $w\ast \Lambda$ for $w\in
\mathcal{W}_{m+n}(\Lambda)$ by definition. It would be interesting
to compare these two formulas.

Note that we can observe  a similarity between ${\rm
ch}L_{m+n}(\Lambda(\lambda))$ and ${\rm
ch}L_{m|n}(\Lambda(\lambda))$ ($\lambda\in\mathbb{Z}_+^d$) from
Theorem \ref{Weyl for hook Schur} and Theorem \ref{Weyl for
unitary}. Recently, a more direct and deep connection between the
Grothendieck groups of two module categories of $\frak{gl}_{m+n}$
and $\frak{gl}_{m|n}$ has been found in \cite{CWZ}, which explains
this similarity of characters.}
\end{rem}

\begin{ex}{\rm Suppose that $\lambda={\bf 0}_d$. Then $\Lambda({\bf 0}_d)=(-d{\bf 1}_m,{\bf
0}_n)$ and each $w\in \mathcal{W}_{m+n}(-d{\bf 1}_m,{\bf 0}_n)$ of
minimal length corresponds to a unique partition $\mu\subset (m^n)$
with $\ell(w)=|\mu|$ (see Remark \ref{comments on cosets} (1)).
Hence, we have $w\ast(-d{\bf 1}_m,{\bf 0}_n)=(-(d^m)-\mu'^*,\mu)$,
and
\begin{equation*}
{\rm ch}L_{m+n}(-d{\bf 1}_m,{\bf 0}_n)=\sum_{\mu\subset
(m^n)}(-1)^{|\mu|} \frac {s_{\mu'+(d^m)}({\bf
x}_{[-m]}^{-1})s_{\mu}({\bf x}_{[n]})}{\prod_{i,j}(1-x_i^{-1}x_j)},
\end{equation*}
where $i\in [-m]$ and $j\in [n]$. }
\end{ex}

\begin{rem}{\rm  We may apply Theorem \ref{main result} to other highest weight
representations of a Lie (super)algebras $\frak{g}$, whenever we
have a Howe duality of $(\frak{g},\frak{gl}_d)$ (cf.
\cite{CW01,CW03,Fr,H,KacR2}), since the associated irreducible
characters for $\frak{g}$ satisfy a Cauchy type identity of the form
given in Theorem \ref{main result} and hence they are equal to
$S_{\lambda}^{\A/\B}$ under suitable choices of $\A$ and $\B$
\cite{Kwon}. For example,
$S_{\lambda}^{\mathbb{Z}'_{>0}/\mathbb{Z}'_{<0}}$ yields a character
of an integrable highest weight representation of
$\widehat{\frak{gl}}_{\infty}$, which is an irreducible component in
fermionic Fock space representations, and
$S_{\lambda}^{\mathbb{Z}_{>0}/\mathbb{Z}_{<0}}$ a character of a
highest weight representation of $\widehat{\frak{gl}}_{\infty}$ (not
integrable) appearing in bosonic Fock space representations
(cf.\cite{Fr,KacR2}). In both cases, the characters are given as an
alternating sum of generalized Verma modules induced from integrable
highest weight representations of a parabolic subalgebra of
$\widehat{\frak{gl}}_{\infty}$ (see Proposition \ref{Weyl type 2}
and \ref{Weyl type 3}). We also have similar applications to super
cases studied in \cite{CW01,CW03,Kwon2} (see also \cite{Kwon}). }
\end{rem}

\end{document}